\documentclass[11pt]{amsart}
\usepackage[latin1]{inputenc}
\usepackage{amsmath,amsthm, amscd, amssymb, amsfonts}
\usepackage[dvips]{graphicx}
\numberwithin{equation}{section}
\newcommand{\s}{\mathbb{S}}
\newcommand{\End}{\mbox{\rm End\,}}
\theoremstyle{plain}

\newtheorem{teor}{Theorem}[section]
\newtheorem{corol}[teor]{Corollary}
\newtheorem{prop}[teor]{Proposition}
\newtheorem{lem}[teor]{Lemma}

\theoremstyle{definition}
\newtheorem{defin}[teor]{Definition}

\newtheorem{ejem}[teor]{Example}

\theoremstyle{remark}
\newtheorem{obs}[teor]{Remark}

\newcommand\id{\operatorname{id}}

\newcommand\co{\operatorname{co}}

\newcommand\ff{\operatorname{ff}}

\newcommand\Hom{\operatorname{Hom}}

\newcommand\Ind{\operatorname{Ind}}

\newcommand\can{\operatorname{can}}
\newcommand\vect{\operatorname{Vec}}
\newcommand\Quot{\operatorname{Quot}}
\newcommand\Sub{\operatorname{Sub}}

\begin{document}

\title[Normal Hopf subalgebras in deformations of finite groups]{Normal Hopf subalgebras in cocycle deformations
of finite groups}
\author{C\' esar Galindo and Sonia Natale}
\address{Facultad de Matem\'atica, Astronom\'\i a y F\'\i sica,
Universidad Nacional de C\'ordoba, CIEM -- CONICET, (5000) Ciudad
Universitaria, C\'ordoba, Argentina}
\email{galindo@mate.uncor.edu, natale@mate.uncor.edu}
\thanks{This work was partially supported by CONICET, Fundaci\' on
Antorchas, Agencia C\'ordoba Ciencia, ANPCyT    and Secyt (UNC)}
\date{November 21, 2007}

\begin{abstract}
Let $G$ be a finite group and let $\pi: G \to G'$ be a surjective
group homomorphism. Consider the cocycle deformation $L =
H^{\sigma}$ of the Hopf algebra $H = k^G$ of $k$-valued linear
functions on $G$, with respect to some convolution invertible
2-cocycle $\sigma$. The (normal) Hopf subalgebra $k^{G'} \subseteq
k^G$ corresponds to a Hopf subalgebra $L' \subseteq L$. Our main
result is an explicit necessary and sufficient condition for the
normality of $L'$ in $L$.
\end{abstract}

\maketitle
\section{Introduction and Main Results}

Let $H$ be Hopf algebra over the field $k$. The category
${}^H\mathcal M$ of its finite dimensional corepresentations is a
tensor category. Tensor categories of this form are characterized,
via tannakian reconstruction arguments, as those possessing a
fiber functor with values in the category $\vect$ of finite
dimensional vector spaces over $k$. The corepresentation
categories of the finite dimensional Hopf algebras $H$ and $L$ are
equivalent if and only if $L = H^{\sigma}$ is a cocycle
deformation of $H$ in the sense of \cite{doi}, that is,  $L = H$
as a coalgebra, and the multiplication in $L$ is given by
$$h._{\sigma}g = \sigma(h_1, g_1)h_2g_2\sigma^{-1}(h_3, g_3),$$ for
some convolution invertible normalized 2-cocycle $\sigma: H
\otimes H \to k$. Equivalently, the dual Hopf algebra $L^*$ is
obtained from $H^*$ via a twisting of the comultiplication through
$\sigma \in H^* \otimes H^*$.

\medbreak A Hopf subalgebra of a finite dimensional Hopf algebra
$H$ is called \emph{normal} if it is invariant under the (left or
right) adjoint action of $H$; $H$ is called \emph{simple} if it
contains no proper normal Hopf subalgebra. For instance, the group
algebra of finite simple group is an example of a (trivial) simple
Hopf algebra.

In a previous paper \cite{GN} we gave a series of examples showing
that the notions of simplicity and (semi)solvability of a
(semisimple) Hopf algebra are \emph{not} twist invariants; that
is, they are not categorical notions. The examples in \textit{loc.
cit.} were obtained as twisting deformations of finite groups;
namely we proved that there exist finite groups which are not
simple, and even supersolvable in certain cases, that admit
twisting deformations which are simple as Hopf algebras.

In this paper we determine necessary and sufficient conditions, in
group-theoretical terms, for a quotient of a twisting deformation
of a finite group to be normal. Any twisting deformation (of the
comultiplication) preserves quotient Hopf algebras. We consider
here the problem of a twisting preserving normality of such a
quotient. We use for this the classification of twisting in group
algebras due to Movshev and Davydov \cite{Mv, Davydov}  and the
results of Schauenburg on Galois correspondences \cite{scha-90}.

\medbreak Another way of regarding twisting deformations is
through the language of Hopf Galois extensions. If $A$ is an
$H$-Galois extension of $k$, also called an \emph{$H$-Galois
object}, then there exists a Hopf algebra $L$, called the left
Galois Hopf algebra of $A$, such that $A$ is an $(L, H)$-bigalois
object and moreover the functor $\mathcal F_A: {}^L\mathcal M \to
{}^H\mathcal M$, $\mathcal F_A(V)= A\square_H V$, is a tensor
equivalence. Conversely, if ${}^L\mathcal M$ and ${}^H\mathcal M$
are tensor equivalent, then there exists an $(L, H)$-bigalois
object $A$ such that the equivalence is essentially $\mathcal F_A$
\cite{galois survey}.

In this context, the Hopf algebra $L$ is a cocycle deformation of
$H$, $L = H^{\sigma}$, and the $H$-Galois object can be identified
with $A = H^{(\sigma)}$, where $H^{(\sigma)} = H$ as right
$H$-comodules with multiplication $h.g = \sigma(h_1, g_1)h_2g_2$,
$h, g \in H$.

\medbreak We shall work over an algebraically closed field $k$ of
characteristic zero, although most results hold as well under
weaker assumptions on the characteristic.

\medbreak Let $G$ be a finite group and let $\pi: G \to G'$ be an
epimorphism of groups. This corresponds to an epimorphism of Hopf
algebras $kG \to kG'$, and in turn to  Hopf algebra inclusion $H'
= k^{G'} \subseteq H = k^G$, which is clearly normal because $k^G$
is commutative. Let $F \subseteq G$ be the kernel of $\pi$; so
that $F$ is a normal subgroup of $G$.

Galois objects for $k^G$ are classified by (conjugacy classes of)
pairs $(S, \alpha)$, where $S \subseteq G$ is a subgroup and
$\alpha \in H^2(S, k^*)$ is a non-degenerate 2-cocycle \cite{Mv,
Davydov}. Let $A = A(G, S, \alpha)$ be the $k^G$-Galois object
associated to this data and let $L = L(A, k^G)$ be the left Galois
Hopf algebra. So that $L^*$ is a twisting deformation of $kG$.
There is a bijective correspondence between Hopf subalgebras $H'
\subseteq k^G$ and Hopf subalgebras of $L$. We consider the Hopf
subalgebra $L' \subseteq L$ corresponding to $H' \subseteq H$.
Then $A' = {}^{\co L/L(L')^+}A \subseteq A$ is a $k^G$-subcomodule
algebra. By results of Schauenburg, $L'$ is normal in $L$ if and
only if $A'$ is stable under the Miyashita-Ulbrich action of $k^G$
on $A$ \cite{scha-90}.

\medbreak For $g \in G$, and $F\subseteq G$ a subgroup, we denote
by $C_F(g) = C_G(g)\cap F$ the centralizer of $g$ in $F$. Suppose
$F\trianglelefteq G$. An element $s \in S$ will be called
$(\alpha,F)$-\emph{regular} if $\alpha(s, t)=\alpha(t, s)$, for
all $t \in C_S(s)\cap F$. The following theorem is proved in
Section \ref{condicion}.

\begin{teor}\label{main} The Hopf subalgebra $L' \subseteq L$ is normal if
and only if $F \subseteq C_G(s)$, for every $(\alpha, F)$-regular
element $s \in S$. \end{teor}

The proof of the theorem reduces to determining under what
conditions the subalgebra of invariants $A^F$ is stable under the
Miyashita-Ulbrich action of $k^G$, where $F$ is the kernel of the
epimorphism $G \to G'$. We find these conditions, more generally,
for the case where $F$ is any (not necessarily normal) subgroup of
$G$; see Theorem \ref{condicion ser submod Miyashita}.

Using Theorem \ref{main} we prove generalizations of some of the
results in \cite{GN}; see Theorem \ref{normalidad para indice
primo} and Corollary \ref{simetrico}. We also give an alternative
proof of the simplicity of a family of twisting deformations of
certain supersolvable groups in Theorem \ref{supersoluble}.

\medbreak The paper is organized as follows: in Section
\ref{prels} we give a brief account of results on Hopf Galois
extensions that will be needed in the sequel. In particular, we
describe in Subsection \ref{correspondence} the Galois
correspondence of Schauenburg and the conditions under which
normal Hopf subalgebras are preserved under this correspondence.
In Section \ref{groups} we review the classification of Hopf
Galois objects for finite groups due to Davydov and Movshev
\cite{Davydov, Mv}; in this context, we give a description of the
Miyashita-Ulbrich action. In Section \ref{invariantes} we
characterize the subalgebra of invariants in a Galois object $A$
under the action of a subgroup. This is applied in Section
\ref{condicion} to give a necessary and sufficient condition for a
quotient of a twisting of a finite group to be normal. In Section
\ref{ejemplos} we apply the previous results to some families of
examples.

\subsection*{Acknowledgement} The authors thank M. Mombelli for
his interest and valuable remarks on this paper.

\section{Preliminaries on Hopf Galois extensions}\label{prels}

In this section we review some results on Hopf Galois extensions
that we shall need later. We refer the reader to \cite{galois
survey} for a detailed exposition on the subject.

\begin{defin}\label{defin hopf galois}
Let $H$ be a Hopf algebra. Let also $A$ be a right $H$-comodule
algebra with structure map $\rho: A\to A\otimes H$, $\rho(a) =
a_{(0)}\otimes a_{(1)}$, and let $B:= A^{\co H}$.  The extension
$B\subseteq A$ is called a right \emph{Hopf Galois} extension, or
a right $H$-\emph{Galois}  extension if the canonical map $$\can:
A\otimes_BA\to A\otimes H, \  \  x\otimes y\mapsto xy_{(0)}\otimes
y_{(1)},$$is bijective.

A right $H$-Galois extension of the base field $k$ will be called
a right $H$-\emph{Galois object}. Left $H$-Galois extensions and
left $H$-Galois objects are defined similarly.
\end{defin}

\begin{defin}
An $(L, H)$-bigalois object is an $(L, H)$-bicomodule algebra $A$
which is simultaneously a left $L$-Galois object and right
$H$-Galois object.
\end{defin}

\begin{ejem}\label{ejem1}
Let the Hopf algebra $H$ coact on itself through the
comultiplication $\Delta$. Then the canonical map  $ H\otimes H\to
H\otimes H,  \ x\otimes y\mapsto xy_{(1)} \otimes y_{(2)}$ is
bijective with inverse $a\otimes h \mapsto aS(h_{(1)})\otimes
h_{(2)}$. Since $\Delta$ is coassociative, $H$ is an $(H,
H)$-bigalois object.

More generally, let $\pi:H\to Q$ be an epimorphism of finite
dimensional Hopf algebras. Then $H$ is a right $Q$-comodule
algebra with structure map $\rho=(\id\otimes\pi)\circ \Delta$. By
\cite[Theorem 8.2.4 and Proposition 8.4.4]{Mont} $H^{\co
Q}\subseteq H$ is a right $Q$-Galois extension. \end{ejem}

\emph{Assume from now on that $H$ is  finite dimensional.} We have
a characterization of  $H^*$-Galois objects $A$ in terms of the
natural $H$-action over $A$:  $A$ is a right $H^*$-Galois object
if and only if $A$ is finite dimensional and the map $\pi: A\#H\to
\End A$, $\pi(a \# h)(b) = a(h.b)$, is an isomorphism.

\begin{prop}\label{semisimple galois}
Every right $H$-Galois object is isomorphic to a crossed product
$H^{(\sigma)} = k\#_{\sigma}H$, where $\sigma:H\otimes H\to k$ is
an invertible 2-cocycle. Moreover, if $H$ is cosemisimple, then
$H^{(\sigma)}$ is semisimple.
\end{prop}

On the other hand, $H^{(\sigma)}$ is a right $H$-Galois object,
for all such $\sigma$.

\begin{proof} Let $A$ be a right $H$-Galois object.
By  \cite[Proposition 8.3.6]{Mont}, $A\simeq k\#_{\sigma}H$, where
$\sigma: H\otimes H\to k$ is an invertible 2-cocycle, since $A \#
H^* \simeq \End A$ is simple artinian. Now by \cite[Theorem
7.4.2]{Mont} if $H$ is cosemisimple and finite dimensional,
$A\simeq k\#_{\sigma}H$ is semisimple.
\end{proof}

For any right $H$-Galois object $A$ there is an associated Hopf
algebra $L = L(A, H)$, called the \emph{left Galois} Hopf algebra,
such that $A$ is in a natural way an $(L, H)$-bigalois object.

In the context of Proposition \ref{semisimple galois}, this Hopf
algebra is the \emph{cocycle deformation} $H^{\sigma}$ of $H$,
obtained from $H$ keeping the coalgebra structure unchanged and
deforming the multiplication via $$h._{\sigma}g = \sigma(h_1,
g_1)h_2g_2\sigma^{-1}(h_3, g_3).$$See \cite[Proposition
3.1.6]{galois survey}.

In the dual situation, we have that $L^* = (H^*)^{\sigma}$ is
obtained from $H^*$ through a twisting of the comultiplication:
$\Delta_{\sigma}(x) = \sigma \Delta(x)\sigma^{-1}$, $x \in L^*$.

 \subsection{Miyashita-Ulbrich action}\label{miy-ul} Let $B \subseteq
A$ be an $H$-Galois extension. The Miyashita-Ulbrich action of $H$
on the centralizer $A^B$ of $B$ in $A$, makes $A^B$ into a
commutative algebra in $\mathcal{YD}_H^H$.

\begin{defin}
Let $H$ be a Hopf algebra, and $A$ an $H$-Galois extension of $B$.
The Miyashita-Ulbrich action of $H$ on $A^B$ is defined by
$x\leftharpoonup h = h^{[1]}xh^{[2]}$, $x\in A^B$, $h\in H$, where
$\can(h^{[1]}\otimes h^{[2]}) = 1\otimes h$.
\end{defin}

The Miyashita-Ulbrich action of $H$ on the $H$-Galois object $A$
is characterized as the unique map $A\otimes H\to A$, $x\otimes
h\mapsto x\leftharpoonup h$, such that $xy=y_{(0)}(x\leftharpoonup
y_{(1)})$, for all $x, y \in A$.

\begin{ejem} Recall that a Hopf subalgebra $H' \subseteq H$ is called
\emph{normal} if it is stable under the (left or right) adjoint
action of $H$.

Consider the $H$-Galois object $A=H$ as in Example \ref{ejem1}. In
this case, the Miyashita-Ulbrich action coincides with the right
adjoint action of $H$ on itself: $x.h = \mathcal
S(h_{(1)})xh_{(2)}$, $x, h \in H$.

Therefore $H' \subseteq H$ is a normal Hopf subalgebra if and only
if it is stable under the Miyashita-Ulbrich action. In the next
subsection we recall some results of Schauenburg that we shall use
later and which generalize this situation. \end{ejem}

\subsection{Galois correspondences}\label{correspondence}

Let $L$, $H$ be finite dimensional Hopf algebras. Let $A$ be a
nonzero $(L, H)$-bigalois object.

\medbreak Let $\Quot(L)$  and $\Sub^H_{\ff}(A)$ be, respectively,
the set of coideal left ideals $I \subseteq L$, and the set of
$H$-comodule subalgebras $B \subseteq A$ such that $A$ is left
faithfully flat over $B$. In view of results of Schauenburg
\cite[Theorem 3.6]{scha-90} there are mutually inverse bijections
\begin{equation}\label{gal-corr}\Quot(L) \leftrightarrows
\Sub^H_{\ff}(A),\end{equation} defined as follows: $I \in
\Quot(L)$ is assigned to ${}^{\co L/I}A \subseteq A$, and $B
\subseteq A$ goes to the coideal left ideal $\mathcal J(B)$ with
$L/\mathcal J(B) : = (A \otimes_BA)^{\co H}$. We note that
faithful flatness assumptions in \cite{scha-90} are guaranteed in
the finite dimensional context by the results of Skryabin
\cite{skryabin}.

\medbreak By \cite[Corollary 3.13]{scha-90} there is a bijective
correspondence between Hopf subalgebras of $H$ and $L$. Let $H'
\subseteq H$ be a Hopf subalgebra. Then $A(H') = A \square_H H'
\subseteq A \square_H H \simeq A$ is a right $H'$-Galois object.
Let $L' = L(A(H'), H')$. Then  $L' \subseteq L$ is a Hopf
subalgebra in a unique way such that $A(H') \subseteq A$ is an
$(L, H)$-bicomodule subalgebra. In this case we have that $A(H')$
is an $(L', H')$-bigalois object, and $A(H') = {}^{\co \overline
L}A$, $L' = {}^{\co \overline L}L$, $\overline L = L / L(L')^+$.

\medbreak Suppose $H' \subseteq H$ is a Hopf subalgebra. The
corresponding Hopf subalgebra $L' \subseteq L$ gives rise in turn
to a coideal left ideal $I = L(L')^+ \in \Quot(L)$, which
corresponds, under \eqref{gal-corr}, to the $H$-comodule
subalgebra ${}^{\co L/I}A = A(H')$.

By \cite[Theorem 3.8 (i)]{scha-90}, $I = L(L')^+ \in \Quot(L)$ is
a \emph{Hopf ideal}, that is, $L' \subseteq L$ is a \emph{normal}
Hopf subalgebra, if and only if $A(H') \subseteq A$ is an
$H$-submodule under the Miyashita-Ulbrich action. See also
\cite[pp. 593]{scha}.

\medbreak The following lemma describes the subalgebra $A(H')$.

Consider the quotient left $H$-module coalgebra $p: H \to H''$,
where $H'' = H/H(H')^+$. The coalgebra $H''$ coacts on $A$ via
$(\id \otimes p)\rho: A \to A \otimes H''$.

\begin{lem}\label{desc-a'} We have $A(H') =
\rho^{-1}(A\otimes H') = A^{\co H''}$. \end{lem}

Here we identify $A(H')$ with a subalgebra of $A$ through the
canonical isomorphism $\rho: A \to A\square_HH$.

\begin{proof} It is clear that under the isomorphism $\rho: A \to A\square_HH$,
$A\square_H H'$ is identified with $\rho^{-1}(A\otimes
H')\subseteq A$. The inclusion $A(H') \subseteq A^{\co H''}$
follows since $p\vert_{H'} = \epsilon$. The other inclusion
follows from the fact that $H' = H^{\co H''}$.  This finishes the
proof of the lemma.
\end{proof}

\begin{obs} The coaction $\rho: A \to A \otimes H$ gives rise to a left
action $H^* \otimes A \to A$, such that $f.a = f(a_{(1)})a_{(0)}$.

Consider the right coideal subalgebra $R = (H'')^* \subseteq H^*$.
Then we have $A(H') = A^{\co H''} = A^R$, where $A^R$ is the
subalgebra of invariants in $A$ under the action of $R$.\end{obs}

\section{Galois extensions of finite groups}\label{groups}

Let $G$ be a finite group, and let $A$ be a $G$-module algebra,
with the $G$-action $kG \otimes A \to A$, $g \otimes a \mapsto
g.a$. Equivalently, $A$ is a $k^G$-comodule algebra with respect
to the coaction $\rho: A \to A \otimes k^G$ defined by $$\rho(a) =
\sum_{g \in G} g.a \otimes e_g.$$ Here, and elsewhere in this
paper, $(e_g)_{g \in G}$ denotes the basis of $k^G$ consisting of
the canonical idempotents $e_g(h) = \delta_{g, h}$, $g, h \in G$.

\medbreak Let $S$ be a finite group and let $\alpha \in Z^2(S,
k^*)$ be a 2-cocycle.
 For each $s \in S$, we shall use the notation $x_s \in
k_{\alpha}S$ to indicate the corresponding element in the twisted
group algebra $k_{\alpha}S$. Thus $(x_s)_{s \in S}$ is a basis of
$k_{\alpha}S$ and, in this basis, $x_sx_t = \alpha(s, t)x_{st}$.

\smallbreak Recall that an element $s \in S$ is called
$\alpha$-regular if $\alpha(s, t) = \alpha(t, s)$, for all $t \in
C_S(s)$. This definition depends only on the class of $s$ under
conjugation. The 2-cocycle $\alpha$ is called
\emph{non-degenerate} if and only if $\{1\}$ is the only
$\alpha$-regular class in $S$.

Suppose $\alpha$ is a non-degenerate 2-cocycle on $S$ and the
characteristic of $k$ is relatively prime to $|S|$. Then the
twisted group algebra $k_{\alpha}S$ is naturally a $k^S$-Galois
object with left $S$-action given by $s \triangleright x_t =
x_sx_tx_s^{-1}$. This Galois object depends only on the cohomology
class of $\alpha$ in $H^2(S, k^*)$ \cite[Proposition
3.5]{Davydov}. In this case, the Miyashita-Ulbrich action
$k_{\alpha}S \otimes k^S \to k_{\alpha}S$ corresponds to the
standard $S$-grading of $k_{\alpha}S$, with respect to which $x_t$
is homogeneous of degree $t \in S$.

\medbreak We note for future use:

\begin{lem}\label{normalizacion}\cite[Lemma 3.6]{karpi}\label{cociclo especial}
Let $G$ be a finite group.  Then any $\alpha \in Z^2(G,k^*)$ is
cohomologous to $\beta \in Z^2(G,k^*)$ such that
$\beta(x,x^{-1})=1$, for all $x\in G$. \qed
\end{lem}

We may and shall assume from now on that the non-degenerate
2-cocycle $\alpha\in Z^2(S,k^*)$ satisfies $\alpha(s,s^{-1})=1$,
for all $s\in S$.

This assumption implies that, in the twisted group algebra
$k_{\alpha}S$, we have $x_t^{-1} = x_{t^{-1}}$, for all $t \in S$.

\begin{teor}\label{davydov}\cite{Davydov, Mv}.
Isomorphism classes of $k^G$-Galois objects are in one-to-one
correspondence  with conjugacy classes of  pairs $(S,\alpha)$,
where $S$ is a subgroup of $G$ and $\alpha\in Z^2(S,k^*)$ is a
non-degenerate 2-cocycle.

The $k^G$-Galois object corresponding to the pair $(S,\alpha)$ can
be constructed as the algebra of $S$-invariant functions
$$A(G,S,\alpha)=\{r:G\to k_\alpha S:\  \  f(sg)=s\rhd f(g)\},$$
where $k_\alpha S$ is the twisted group algebra, with action
$s\rhd x_t = x_sx_tx_s^{-1}$, $s, t \in S$, and the $G$-action on
$A$ is $(g.f)(h) = f(hg)$.\qed
\end{teor}

The following proposition gives an alternative description of the
Galois object $A(G, S, \alpha)$ that will be convenient to our
purposes.

\begin{prop}\label{identif}
There is an isomorphism of $G$-algebras $$A(G, S, \alpha) \simeq
\Ind_S^G k_{\alpha}S = kG\otimes_{kS}k_{\alpha}S,$$ where the left
$G$-action and the multiplication are given, respectively, by
\begin{align} \label{f1}g'. g \otimes x & = g'g \otimes x, \\
\label{f2}(g \otimes x) (h \otimes y) & = \begin{cases}0, \qquad
\qquad \qquad \; \, h^{-1}g \notin S, \\ h \otimes (h^{-1}g . x)y,
\quad h^{-1}g \in S.
\end{cases}
\end{align}
\end{prop}

\begin{proof} Let $A = A(G,S,\alpha)$. By definition, $A = \Hom_{kS}(kG,k_{\alpha}S)$.
Using the version of Frobenius reciprocity in \cite[Proposition
(10.21)]{CR}, we have isomorphisms $\Hom_{kS}(kG, k_{\alpha}S)
\simeq \Hom_{kG}(kG, kG \otimes_{kS}k_{\alpha}S) \simeq kG
\otimes_{kS}k_{\alpha}S$.

Write $G = \cup_{i = 1}^r g_iS$, where $g_i$, $1 \leq i \leq r$,
is a set of representatives of the left cosets of $S$ in $G$. An
isomorphism $\phi: A(G, S, \alpha) \to kG\otimes_{kS}k_{\alpha}S$
is defined explicitly by
$$\phi(f) = \sum_i g_i \otimes f(g_i^{-1}).$$ Its inverse $\psi: kG\otimes_{kS}k_{\alpha}S \to A(G, S,
\alpha)$ is determined by $$\psi(g \otimes x)(b) = \lambda(bg)x,$$
where $\lambda: G \to kS$ is the map $\lambda(g) = g$, if $g \in
S$, and $\lambda(g) = 0$, otherwise. Of course, these isomorphisms
do not depend on the choice of left coset representatives.
Moreover, they transform the $G$-algebra structure on $A$ into the
$G$-algebra structure on $kG\otimes_{kS}k_{\alpha}S$ given by
formulas \eqref{f1} and \eqref{f2}. This proves the proposition.
\end{proof}

From now on we shall use the identification $A = A(G, S, \alpha)
\simeq kG\otimes_{kS}k_{\alpha}S$ given in Proposition
\ref{identif}.

\medbreak We next describe the Miyashita-Ulbrich action
$\leftharpoonup: A \otimes k^G \to A$ in terms of this
identification. First note that, in our context, this action
corresponds to a $G$-grading on $A$: $A = \oplus_{g\in G}A_g$,
such that $A_g = A \leftharpoonup e_g$, $g \in G$. We shall call
it the Miyashita-Ulbrich grading.

\begin{lem}\label{mu-grading} The Miyashita-Ulbrich grading on $A$ is determined by $$|g \otimes x_s| : = g s
g^{-1},$$ for all $g \in G$, $s \in S$. \end{lem}

\begin{proof} Recall that the Miyashita-Ulbrich action of a Hopf algebra $H$ over
the $H$-Galois object $A$ is characterized as the unique map
$A\otimes H\to A$, $x\otimes h\mapsto x\leftharpoonup h$, such
that $xy=y_{(0)}(x\leftharpoonup y_{(1)})$, for all $x, y \in A$.
Thus, for $H = k^G$, we have $A_g = \{a\in A| \, ax = (g\cdot x)a,
\; \text{for all } x\in A\}$.

Let $g, h \in G$, $s, t \in S$. We have
\begin{align*}\left((gsg^{-1}). (h \otimes x_t) \right) (g \otimes
x_s) & = (gsg^{-1}h \otimes x_t)(g \otimes x_s) \\ & =
\begin{cases}0, \qquad \qquad \qquad \; \, sg^{-1}h \notin S, \\ g
\otimes (sg^{-1}h . x_t)x_s, \quad sg^{-1}h \in S.
\end{cases}\end{align*}
Clearly, $sg^{-1}h \in S$ if and only if $h^{-1}g \in S$ and, in
this case, $$g \otimes (sg^{-1}h . x_t)x_s = g \otimes x_s
(g^{-1}h . x_t) = h \otimes (h^{-1}g . x_s) x_t,$$ the first
equality because $x_s$ is homogeneous of degree $s$ with respect
to the Miyashita-Ulbrich $S$-grading in $k_{\alpha}S$. Hence
$$\left((gsg^{-1}). (h \otimes x_t) \right) (g \otimes x_s) = (g
\otimes x_s) (h \otimes x_t).$$ Since this holds for all $g, h \in
G$, $s, t \in S$, then $g \otimes x_s$ is homogeneous of degree
$gsg^{-1}$, as claimed. This finishes the proof of the lemma.
\end{proof}

\section{The subalgebra $A^F$  of $F$-invariants for a subgroup $F$}\label{invariantes}

Following \cite[Chapter 7, Section 9]{karpi} we shall consider
monomial representations in order to describe a basis of the
subalgebra of invariants $A^F$.

\medbreak Let $F$ be a finite group. A \textit{monomial space}
over $k$ consists of the data $(V, X, (V_x)_{x\in X})$, where $V$
is a vector space, $X$ is finite set, and $(V_x)_{x\in X}$ is a
family of one-dimensional subspaces of $V$ indexed by $X$, such
that $V = \oplus_{x\in X}V_x$.

By a \emph{monomial representation} of a group $F$ on $(V,X,
(V_x)_{x\in X})$ we mean a group homomorphism $$\Gamma: F\to
GL(V),$$such that for each $\sigma\in F$, $\Gamma(\sigma)$
permutes the $V_x$'s; hence the action on $V$ induces an action by
permutations of $F$ on $X$.

\medbreak For each $x\in X$, we shall denote $F(x)$ the stabilizer
of $x$. We say that an element $x\in X$ is \emph{regular} under
the monomial action of $F$ if, for all $\sigma\in F(x)$,
$\Gamma(\sigma)$ is the identity map on $V_x$.

\begin{lem}\cite[Lemma 9.1]{karpi}. \label{lema monomial}
Let $(V,X, (V_x)_{x \in X})$ be a monomial space, and let
$\Gamma:F\to GL(V)$ a monomial representation of $F$ on $(V,
X,(V_x)_{x \in X})$.
\begin{enumerate}
  \item If $x\in X$ is a regular element under the monomial action of $F$,
  then so are all elements in the $F$-orbit of $x$.
  \item Let $T$ be a set of representatives of the regular orbits of
  $X$. For each $t\in T$, let $w_t$ be a nonzero vector of $V_t$. Set
  $$v_t=\sum_{\sigma\in Y_t}\Gamma(\sigma)w_t,$$ where $Y_t$ is a set of left coset representatives for $F(t)$ in $F$.
  Then $(v_t)_{t\in T }$ is a basis of the subspace $V^F$ of $F$-invariants  in $V$.
\qed \end{enumerate}
\end{lem}

\begin{obs}\label{linind} The dimension of $V^F$ is equal to the number of regular $F$-orbits
under the monomial action of $F$ on $X$.

Note in addition that the set $\{\Gamma(\sigma)w_t\}_{t\in T,
\sigma \in Y_t}$ is linearly independent in $V$. Indeed,
$\Gamma(\sigma)w_t$ belongs to $V_{\sigma.t}$. Because elements in
$T$ are not conjugated under the action of $F$, and $\sigma$ runs
over a set of left coset representatives of $F(t)$, all these
elements have pairwise distinct degree of homogeneity and are thus
linearly independent. \end{obs}

Let $I=\{g_1, g_2, \ldots, g_n\}$ be a set of representatives of
the right  cosets of $S$ in $G$,  so that
$\{g_i\otimes_{kS}x_s\}_{s\in S,i\in I}$ is a basis of
$kG\otimes_{kS}k_{\alpha}S$. It is clear that $(A, I\times S, k
(g\otimes x_s))$ is a monomial space.

\medbreak The set $I$  carries a natural permutation action of $G$
given by $g.g_i = g_j$, if and only if
$$gg_i = g_jt, \qquad  t \in S.$$ The action of $G$ on $A$ in
terms of the basis $\{g_i\otimes_{kS}x_s\}_{s\in S,i\in I}$, has
the form
$$g. (g_i\otimes_{kS}x_s)= gg_i\otimes_{kS}x_s=
(g. g_i)t \otimes_{kS}x_s=
\alpha(t,s)\alpha(ts,t^{-1})(g.g_i)\otimes_{kS}x_{tst^{-1}}.$$
Then, for each $g\in G$, the action of $g$ permutes the spaces
$k(g_i\otimes_{kS}x_s)$; in other words, $(A, I\times S, k
(g_i\otimes x_s))$ is a monomial representation of $G$.

The induced action of $G$ on $I\times S$ is given by $$g.(g_i, s)
= (g_j, tst^{-1}) = (g.g_i, tst^{-1}),$$ with $g_j\in I$, $t\in
S$, as above.

\medbreak \emph{Let $F\subseteq G$ be a subgroup. In what follows
we shall consider the monomial representation of $F$ on $A$
obtained by restriction.}

\medbreak We can regard $S$ as a left $S$-set and $G$ as a right
$S$-set, with actions given by conjugation and right
multiplication, respectively; that is,  $s . t= sts^{-1}$, $g . s
= gs$ for $s,t\in S,\  \  g\in G$. Then we can consider the
quotient set $G\times_S S$ under the equivalence relation: $(g,
s)\sim (g',s')$, if and only if, there exist $t \in S$ such that
$(g', s')=(g. t^{-1}, t. s) = (g t^{-1}, t st^{-1})$.

The set $G\times_S S$ is equipped with a left $G$-action given by
left multiplication on the first component.

\begin{lem}
The map $I\times S \to G\times_S S$ that carries the element
$(g_i,s)$ to the class of the pair $(g_i, s)$, is an isomorphism
of $G$-sets.
\end{lem}

\begin{proof} For every $g \in G$, write $g = g_is'$, uniquely, for some $g_i \in I$, $s' \in S$.
Map the class of the pair $(g, s)$ to the pair $(g_i,
s's{s'}^{-1}) \in I \times S$. This gives a well defined inverse
of the map in the lemma. It is clear that both maps are maps of
$G$-sets.
\end{proof}

We can identify the monomial representation  $(A, I\times S, k
(g_i\otimes x_s))$ of $F$  with $(A, G\times_S S, k (g\otimes
x_s))$.

\begin{defin}
We shall say that a class $(g,s)\in G\times_S S$ is
$(\alpha,F)$-regular if $\alpha(s, t)=\alpha(t, s)$, for all $t
\in C_S(s)\cap g^{-1}Fg$.

Suppose that $F$ is normal in $G$. We shall say that $s \in S$ is
$(\alpha,F)$-regular if $\alpha(s, t)=\alpha(t, s)$, for all $t
\in C_S(s)\cap F$.
\end{defin}

\begin{obs}\label{regularidad} In view of Lemma \ref{lema F-invariant} below, the notion of $(\alpha, F)$-regularity
is well defined, that is, it depends only on the class of
$(g,s)\in G\times_S S$.

Moreover, under the assumption that $F$ is normal in $G$, the
definition of an $(\alpha,F)$-regular element $s \in S$ depends
only on the conjugacy class of $s$ in $S$. In this case, for all
$g\in G$, the class $(g,s)$ is $(\alpha,F)$-regular in $G\times_S
S$ if and only if $s \in S$ is $(\alpha,F)$-regular.
\end{obs}

\begin{lem}\label{lema F-invariant}
The class $(g, s)\in G\times_S S$ is regular under the action of
$F$ if and only if it is $(\alpha,F)$-regular.
\end{lem}

\begin{proof}
The stabilizer of $(g, s)$ in $F$ is the subgroup $F(g, s) =
gC_S(s)g^{-1}\cap F$. Let $r\in gC_S(s)g^{-1}\cap F$, \  $r=
gtg^{-1}, \  t\in C_S(s)$. Then
\begin{align*}
gtg^{-1}\rhd g\otimes_{kS}x_s &= \alpha(t,s)\alpha(ts,t^{-1})g\otimes_{kS}x_s\\
& =\alpha^{-1}(s,t) \alpha(t,s)g\otimes_{kS}x_s\\
& =\alpha^{-1}(s,g^{-1}rg)\alpha(g^{-1}rg,s)g\otimes_{kS}x_s,
\end{align*}
where the second equality  follows because  $\alpha(g,g^{-1})=1$,
for all $g\in G$. Then the class $(g,s)$ is regular under the
action of $F$ if and only if  $\alpha^{-1} (s, g^{-1}rg)$
$\alpha(g^{-1}rg, s) = 1$, for all $r\in g^{-1}C_S(s)g\cap F$,
\textit{i.e.}, if and only if $(g,s)$ is $(\alpha, F)$-regular.
\end{proof}

\begin{prop}\label{bases invariant}
Let $T$ be a set of representatives of the regular $F$-orbits of
$G\times_S S$. For each $(g, t)\in T$, set
$$v_{(g, t)}=\sum_{h\in Y_{(g, t)}} hg \otimes_{kS}x_t, \qquad  (g, t)\in T,$$
where $Y_{(g, t)}$ is a set of left coset
representatives of $gC_S(s)g^{-1}\cap F$ in $F$. Then $(v_{(g,
t)})_{(g, t)\in T}$ is a basis of the subalgebra $A^F$ of
$F$-invariants of $A$.
\end{prop}

\begin{proof} We have that
$(A,G\times_SS,k(g \otimes_{kS}x_s))$ is a monomial representation
of $F$. By Lemma \ref{lema F-invariant},  the stabilizers are the
subgroups $gC_S(s)g^{-1}\cap F$. By Lemma \ref{lema monomial},
$(v_{(g, t)})_{(g, t)\in T }$ is basis of $A^F$.
\end{proof}

\begin{obs}\label{dimension} Since the $k^G$-Galois object
$A(G,S,\alpha)$ is isomorphic to a crossed product over $k^G$, it
is isomorphic to the regular representation $kG$ as a $G$-module.
Therefore $\dim A(G,S,\alpha)^F= [G:F]$.

In view of Proposition \ref{bases invariant}, for all
nondegenerate cocycle $\alpha$, the number of $(\alpha,
F)$-regular orbits in $G \times_SS$ equals $[G: F]$.
\end{obs}

\section{Normal Hopf subalgebras in cocycle deformations}\label{condicion}

Let $G$ be a finite group, let $A = A(G, S, \alpha)$ be a
$k^G$-Galois object, and consider the twisted Hopf algebra
deformation $L = L(A, k^G)$ of $k^G$.

\medbreak Let $\pi: G \to G'$ be an epimorphism of groups. This
corresponds to an epimorphism of Hopf algebras $kG \to kG'$. Let
$F \subseteq G$ be the kernel of $\pi$; so that $F$ is a normal
subgroup of $G$.

Dualizing the epimorphism $\pi: G \to G'$ we obtain a  Hopf
algebra inclusion $H' = k^{G'} \subseteq H = k^G$, which is of
course normal since $k^G$ is commutative.

 We  determine necessary and sufficient
conditions, in group-theoretical terms, for the Hopf subalgebra
$L' \subseteq L$, corresponding to $H' \subseteq H$, to be normal
in $L$.

\medbreak First, for any subgroup $F$ of $G$, we have:

\begin{teor}\label{condicion ser submod Miyashita}
$A^F$ is a $k^G$-submodule of $A$ with respect to the
Miyashita-Ulbrich action if and only if  for all $(\alpha,
F)$-regular $(g, s) \in G\times_S S$, $g^{-1}Fg \subseteq C_G(s)$.
\end{teor}

\begin{proof}
We keep the notations in Proposition \ref{bases invariant}. Note
that the condition $g^{-1}Fg \subseteq C_G(s)$ is equivalent to
$C_F(gsg^{-1}) = F$.

\medbreak Suppose  first that  $C_F(gsg^{-1}) = F$, for all $(g,
s)\in T$. Let $\sigma \in G$. By Lemma \ref{mu-grading}, we have
\begin{align*}
v_{(g, s)}\leftharpoonup e_{\sigma} & =  \sum_{h\in Y_{(g,s)}} hg\otimes_{kS}x_s\leftharpoonup e_{\sigma}\\
                                     & =  \sum_{h\in Y_{(g,s)}}\delta_{\sigma,(hg)s(hg)^{-1}} hg\otimes_{kS}x_s\\
                                     & =  \delta_{\sigma,gsg^{-1}} v_{(g, s)}. \end{align*}
Since this holds for all $\sigma \in G$, then $A^F$ is a
$k^G$-submodule with respect to the Miyashita-Ulbrich action.

\medbreak Conversely, suppose that $A^F$ is a $k^G$-submodule
under the Miyashita-Ulbrich action. Let $\{v_{(g_0,t_0)},
v_{(g_1,t_1)},\ldots, v_{(g_r,t_r)}\}$ be a basis of $A^F$ given
by Proposition \ref{bases invariant}.

Recall that $\{hg_i\otimes_{kS}x_{t_i}\}_{h\in Y_{(g_i,t_i)},
(g_i,t_i)\in T}$ is a set of linearly independent vectors in $A$;
c.f. Remark \ref{linind}.

We have $v_{(g_0,t_0)}= \sum_{i = 1}^k v_{(g_0,t_0)}^{\sigma_i}$,
where $v_{(g_0,t_0)}^{\sigma_i}$ is the homogeneous component of
degree $\sigma_i \in G$ with the respect to the Miyashita-Ulbrich
grading:
$$v_{(g_0,t_0)}^{\sigma_i} = \sum_{h\in Y_{(g_0,t_0)},\,
hg_{0}t_0g_{0}^{-1}h^{-1} = \sigma_i} hg_{0}\otimes_{kS}x_{t_0}.$$
Since, by assumption, $A^F$ is stable under the Miyashita-Ulbrich
action of $k^G$, each one of the homogeneous components
$v_{(g_0,t_0)}^{\sigma_i}$ is in $A^{F}$.

Thus, if $v_{(g_0,t_0)}^{\sigma_1},\ldots,
v_{(g_0,t_0)}^{\sigma_l}$ are the nonzero homogeneous components
of $v_{(g_0,t_0)}$, then
$$\{v_{(g_0,t_0)}^{\sigma_1},\ldots,
v_{(g_0,t_0)}^{\sigma_l},v_{(g_1,t_1)},\ldots, v_{(g_r,t_r)} \}$$
is a set of $l+r$ linearly independent vectors in $A^F$. But
$v_{(g_0,t_0)}$, $v_{(g_1,t_1)}$, $\ldots, v_{(g_r,t_r)}$ form a
basis of $A^F$. Therefore $l = 1$. Hence $v_{(g_0,t_0)}^{\sigma_1}
= v_{(g_0,t_0)}$, and $v_{(g_0,t_0)}^{\sigma_i} = 0$, for all $i
> 1$. This implies that $v_{(g, s)}$
 is homogeneous for all $(\alpha, F)$-regular class $(g, s)$.

\medbreak Let $(g, s)$ be an $(\alpha, F)$-regular class. We have
$v_{(g, s)}\leftharpoonup e_u = v_{(g, s)}$, for some $u \in G$.
Thus
$$v_{(g, s)}   = \sum_{h\in Y_{(g, s)}} hg \otimes_{kS} x_s \leftharpoonup e_{u}
 = \sum_{h\in Y_{(g, s)}}\delta_{u,(hg) s (hg)^{-1}}
 hg\otimes_{kS}x_s.$$
We may assume that $e\in Y_{(g, s)}$, where $e \in F$ is the
identity element. Then, by linear independence of the set
$\{hg\otimes_{kS}x_{t}\}_{h\in Y_{(g, t)}, (g, t)\in T}$, we
obtain $gsg^{-1}= u$, and also $gsg^{-1}= hgsg^{-1}h^{-1}$, for
all $h\in Y_{(g, s)}$. But $(g, s)$ is an arbitrary element in the
$(\alpha, F)$-regular $F$-orbit and $Y_{(g, s)}$ is any set of
left coset representatives for $gC_S(s)g^{-1}\cap F$ in $F$. Then
$gsg^{-1}= hgsg^{-1}h^{-1}$ for all $h\in F$. That is, $F =
C_F(gsg^{-1})$. This finishes the proof of the theorem.
\end{proof}

As an application of Theorem \ref{condicion ser submod Miyashita}
we give a proof of our main result.

\begin{proof}[Proof of Theorem \ref{main}] By Lemma \ref{desc-a'} and the results in
\ref{miy-ul}, it is enough to determine the conditions under which
the subalgebra of invariants $A^F$ is stable under the
Miyashita-Ulbrich action. This is a special case of Theorem
\ref{condicion ser submod Miyashita}, for the case when $F$ is a
normal subgroup of $G$.
\end{proof}

\section{Examples}\label{ejemplos}

In this section $G$ will be a finite group, $S\subseteq G$ will be
a subgroup, and $\alpha\in Z^2(S,k^*)$ a non-degenerate 2-cocycle.
As before, $A = A(G, S, \alpha)$ will denote the corresponding
$k^G$-Galois object.

We shall denote by $J \in kS \otimes kS \subseteq kG \otimes kG$
the twist corresponding to $A$. So that $(kG)^J \simeq L^*$, where
$L = L(A, H)$.

\medbreak As a first application of Theorem \ref{main} we state
the following lemma. Its interpretation in terms of twists is
well-known; see \cite[Lemma 5.4.1]{ssld}.

\begin{lem}\label{normal contiene a S}
Let $F$ be a normal subgroup of $G$. If $S\subseteq F$, then $A^F$
is a $k^G$-submodule of $A$ with respect to the Miyashita-Ulbrich
action.
\end{lem}

In other words, if $J \subseteq kF \otimes kF$, then $(kF)^J
\subseteq (kG)^J$ is a normal Hopf subalgebra.

\begin{proof}
Let $(g,s)$ be an $(\alpha,F)$-regular class. Then $\alpha(s, t) =
\alpha(t, s)$, for all $t \in C_S(s)\cap g^{-1}Fg= C_S(s)$. Since
$\alpha$ is non-degenerate, then $s = 1$. By Theorem
\ref{condicion ser submod Miyashita}, $A^F$ is a $k^G$-submodule
of $A$ with respect to the Miyashita-Ulbrich action.
\end{proof}

The following Theorem gives a generalization of \cite[Theorem
3.3]{GN} to the case where $S$ is not necessarily abelian.

\begin{teor}\label{normalidad para indice primo}
Suppose that $Z(G)=1$. Let $F$ be a normal subgroup of $G$, such
that $[G:F]$ is prime. Then  $A^F$ is a $k^G$-submodule of $A$
with respect to the Miyashita-Ulbrich action if and only if $S
\subseteq F$.
\end{teor}

In other words, the corresponding quotient $(kG)^J \to
k(G/F)^{\overline J}$ is conormal, that is, dual to a normal
inclusion, if and only if $J \in kF \otimes kF$.

\begin{proof}   Suppose $S\subseteq F$. Then, by Lemma \ref{normal
contiene a S}, $A^F$ is a $k^G$-submodule of $A$ with respect to
the Miyashita-Ulbrich action.

Conversely, suppose that  $A^F$ is a $k^G$-submodule of $A$ with
respect to the Miyashita-Ulbrich action. Let $s \in S$ be an
$(\alpha, F)$-regular element. Then $F\subseteq C_G(s)$.

\medbreak If $F\varsubsetneq C_G(s)$, then $C_G(s)=G$, because
$[G: F]$ is prime. Since $Z(G)=1$, we have $s=1$. If $F = C_G(s)$,
then $\alpha(s, t)=\alpha(t, s)$, for all $t \in C_S(s)\cap F=
C_S(s)$. Since $\alpha$ is non-degenerate, then $s=1$ as well.

\medbreak Therefore every $(\alpha, F)$-regular class is of the
form $(g,1)$. Recall that $\dim A^F$ is equal to the number of
$(\alpha,F)$-regular orbits in $G\times_SS$. Note that $(g,1)\sim
(\sigma,1)$ if and only if $g\in \sigma S$. Then $(g,1)$ is
conjugate to $(\sigma,1)$ under the action of $F$ if and only if
$g \in F \sigma S = \sigma FS$. Hence the number of
$(\alpha,F)$-regular orbits is $[G:FS]$.

\medbreak On the other hand, the dimension of $A^F$ is $[G:F]$,
see Remark \ref{dimension}. Thus $[G:FS]=[G:F]$, and therefore
$S\subseteq F$. This finishes the proof of the theorem.
\end{proof}

\subsection{Simple deformations of a family of non-solvable groups}
Let $F$ be a finite simple non-Abelian group, and let $x\in F$ be
an element of prime order $p$. Let $G= F \rtimes\mathbb Z_p$ be
the semidirect product with respect to the action given by
conjugation by $x$.  Note that $Z(G)=1$.

\begin{prop}\label{cor-nonsolv}
Let $G = F \rtimes\mathbb Z_p$, then $(kG)^J$ is a simple Hopf
algebra if and only if $S\nsubseteq F$. \end{prop}

\begin{proof}
The only non-trivial normal subgroup of $G$ is  $F$. By Theorem
\ref{normalidad para indice primo}, the quotient $(kG)^J\to
k\mathbb Z_p$ is conormal if and only if $S\subseteq F$.
\end{proof}

\begin{obs} For the  subgroup $S=\langle x\rangle \rtimes \mathbb Z_p$,
we have   $S\simeq  \mathbb Z_p\times \mathbb Z_p$. In particular,
every $1\neq \alpha \in H^2(S, k^*) \simeq \mathbb Z_p \times
\mathbb Z_p$ is non-degenerate. That is, the Hopf algebra $k^G$
has at least $p-1$ non-trivial Galois objects $A(G, S, \alpha)$
for which $S \nsubseteqq F$. All the corresponding twisting
deformations are simple in view of Proposition \ref{cor-nonsolv}.
\end{obs}

As a consequence we get the following generalization of the
results in \cite[Section 3]{GN} for the symmetric groups $\mathbb
S_n$.

\begin{corol}\label{simetrico} Let $n \geq 5$ and let $G = \mathbb S_n$, $F = \mathbb A_n$ be the symmetric and alternating
groups on $n$ letters, respectively. For every subgroup
$S\subseteq \s_n$ and non-degenerate 2-cocycle $\alpha$ on $S$,
the associated deformation $(k\s_n)^J$ is simple if and only if $S
\nsubseteqq \mathbb A_n$. \qed \end{corol}

\begin{proof} We have $G = \mathbb A_n\rtimes \mathbb Z_2 = \s_n$.
The statement follows from \ref{cor-nonsolv}.
\end{proof}

\subsection{Example of family of super solvable groups}
As a further application of Theorem \ref{main} we give in this
Subsection an alternative proof of the main result in
\cite[Section 4]{GN} (note that we drop here the assumption on the
characteristic of $k$ made in \textit{loc. cit.}). Let $p$, $q$
and $r$ be prime numbers such that $q$ divides $p-1$ and $r-1$.
Let $G_1 = \mathbb Z_p \rtimes \mathbb Z_q$ and $G_2 = \mathbb Z_r
\rtimes \mathbb Z_q$ be the only nonabelian groups of orders $pq$
and $rq$, respectively. Let $G = G_1 \times G_2$ and let $\mathbb
Z_q \times \mathbb Z_q \simeq S \subseteq G$ a subgroup of order
$q^2$. In particular, $G$ is supersolvable and $Z(G) = 1$.

\begin{teor}\label{supersoluble} Let $1\neq \alpha \in H^2( S, k^*)$.
Then  $(kG)^J$ is a simple Hopf algebra.
\end{teor}

\begin{proof}
Note that the groups $G_i$ have only one normal subgroup. The
subgroup $S$ is conjugated to a subgroup of the form $S=S_1\times
S_2$, where $S_i\subseteq G_i$,  $S_i\simeq \mathbb Z_q$, and  if
$s=(s_1,s_2)\in S$,  then $C_G(s)=C_{G_1}(s_1)\times
C_{G_2}(s_2)$.

\medbreak Consider a normal subgroup $N$ of $G$ such that $A^N$ is
a $k^G$-submodule of $A$ with respect to the Miyashita-Ulbrich
action.

\medbreak Let  $s \in S$ be an $(\alpha,N)$-regular element; so
that $N\subseteq C_{G}(s)$ and  $\alpha(s, t)=\alpha(t, s)$ for
all $t \in C_G(s)\cap N$.

\medbreak If $s_1\neq 1$ and $s_2\neq 1$ then $C_{G}(s)= S$. But
$S$ contains no normal subgroups of $G$, hence $s_1=1$ or $s_2=1$.
Suppose that $s_1=1$ and $s_2\neq 1$, then $C_{G}(s)= S_1\times
G_2\supseteq N$, so $\alpha(s, t)=\alpha(t, s)$ for all $t \in S$.
But $\alpha$ is non-degenerate, so $s=1$, which is a
contradiction.

\medbreak Then every $(\alpha, F)$-regular class is of the form
$(g,1)$. The same argument of the proof of Theorem \ref{normalidad
para indice primo} shows that $[G:NS]=[G:N]$, and thus $S\subseteq
N$. This implies that $N=G$. Hence $(kG)^J$ contains no proper
normal Hopf subalgebra. This proves the theorem.
\end{proof}

\bibliographystyle{amsalpha}

\begin{thebibliography}{A}

\bibitem{CR} {\sc C. Curtis} and {\sc I. Reiner}, \emph{Methods of representation theory, I},
Wiley Interscience Publications, New York (1990).

\bibitem{Davydov} {\sc A. A. Davydov}, \emph{Galois Algebras and Monoidal Functors between
Categories of Representations of Finite Groups},  J. Algebra \textbf{244} (2001), 273--301.

\bibitem{doi} {\sc Y. Doi}, \emph{Braided bialgebras and quadratic bialgebras}, Commun. Algebra
\textbf{21} (1993), 1731--1749.

\bibitem{GN} {\sc C. Galindo} and {\sc S. Natale}, \emph{Simple Hopf algebras and deformations of finite groups},
to appear in Math. Res. Lett. Preprint
\texttt{arXiv:math/0608734v2}.

\bibitem{karpi} {\sc G. Karpilovsky}, \emph{Proyective Representation of Finite Groups},
Pure and Applied Mathematics \textbf{94}, Marcel Dekker, New
York-Basel (1985).

\bibitem{ssld} {\sc S. Natale}, \emph{Semisolvability of semisimple Hopf algebras of low dimension},
Memoirs Amer. Math. Soc. \textbf{186} (2007).

\bibitem{Mv} {\sc M. Movshev}, \emph{Twisting in group algebras of finite groups},  Func. Anal.
Appl. {\bf 27} (1994), 240--244.

\bibitem{Mont} {\sc S. Montgomery}, \emph{Hopf Algebras and Their Action on Rings},  CBMS \textbf{82},  Am. Math.
Soc., Providence, Rhode Island (1993).

\bibitem{galois survey} {\sc P. Schauenburg},
\emph{Hopf-Galois and Bi-Galois Extensions}, Fields Institute
Communications \textbf{43} (2004), 469--515.

\bibitem{scha-90} {\sc P. Schauenburg}, \emph{Galois correspondence for Hopf Bigalois extensions},
J. Algebra \textbf{201} (1990), 53--70.

\bibitem{scha} {\sc P. Schauenburg}, \emph{Galois objects over generalized Drinfeld doubles, with an application to
$u_q(sl_2)$}, J. Algebra \textbf{217} (1999), 584--598.

\bibitem{skryabin}  {\sc S. Skryabin},  \emph{ Projectivity and
freeness over comodule algebras}, Trans. Amer. Math. Soc.
\textbf{359}  (2007), 2597--2623.

\end{thebibliography}

\end{document}